\documentclass[a4paper,10pt]{amsart}

\usepackage[utf8]{inputenc}
\usepackage[usenames,dvipsnames]{pstricks}

\makeindex

\newtheorem{thm}{Theorem}[section]
\newtheorem{lemma}[thm]{Lemma}

\theoremstyle{definition}
\newtheorem{dfn}[thm]{Definition}
\newtheorem{rem}[thm]{Remark}
\newtheorem{example}[thm]{Example}

\newcommand{\eps}{\varepsilon}
\newcommand{\R}{\mathbb{R}}
\newcommand{\N}{\mathbb{N}}
\newcommand{\grad}{\mathrm{grad}}
\newcommand{\vol}{\mathrm{vol}}
\newcommand{\area}{\mathrm{area}}
\renewcommand{\div}{\mathrm{div}}
\renewcommand{\phi}{\varphi}

\newcommand{\viol}{WOYMP-violating}

\newcommand{\ep}{\,\,}

\begin{document}

\title{Stochastic completeness and volume growth}

\author{Christian B\"ar}
\date{\today}
\address{Universit\"at Potsdam, Inst.\ f.\ Mathematik, Am Neuen Palais 10, 14469 Potsdam, Germany}
\email{baer@math.uni-potsdam.de}
\thanks{This work has been supported by CNPq-CAPES and by Sonderforschungsbereich 647 funded by Deutsche Forschungsgemeinschaft}

\author{G.~Pacelli Bessa}
\address{UFC, Departamento de Matematica, Bloco 914, Campus do Pici, 60455-760 Fortaleza, Ceara, Brazil}
\email{bessa@mat.ufc.br}

\begin{abstract}
It has been suggested in 1999 that a certain volume growth condition for geodesically complete Riemannian manifolds might imply that the manifold is stochastically complete.
This is motivated by a large class of examples and by a known analogous criterion for recurrence of Brownian motion.
We show that the suggested implication is not true in general.
We also give counter-examples to a converse implication.
\end{abstract}

\maketitle \thispagestyle{empty}

\section{Introduction}

Let $M$ be a geodesically complete connected Riemannian manifold.
The Laplace-Beltrami operator $\Delta = \div\circ\grad$ acting on $C^\infty_c(M)$, the space of smooth functions with compact support, is symmetric with respect to the $L^2$-scalar product.
It is well-known that $\Delta$ is essentially self-adjoint in the Hilbert space $L^2(M)$, see e.g.\ \cite[Thm.~5.2.3]{D}.
We denote its unique extension again by $\Delta$.
By functional calculus we can form $e^{t\Delta}$, a bounded self-adjoint operator on $L^2(M)$ for $t\geq 0$.
For any $u_0 \in L^2(M)$ the function $u(x,t) := (e^{t\Delta}u_0)(x)$ solves the \emph{heat equation}
$$
\frac{\partial u}{\partial t} \ep =\ep \Delta u,
$$
$$
u(\cdot,0) \ep= \ep u_0 .
$$
Elliptic regularity theory shows that $e^{t\Delta}$ is smoothing for $t>0$.
Hence there exists $p\in C^\infty((0,\infty)\times M \times M)$ such that 
$$
e^{t\Delta} v (x)\ep =\ep \int_M p(t,x,y)\, v(y)\, dy .
$$
The function $p$ is called the \emph{heat kernel} of $M$.
It has the following properties:
\begin{eqnarray}
p(t,x,y) &>& 0, \nonumber\\
\frac{\partial p}{\partial t} &=& \Delta_x p, \nonumber\\
p(t,x,y) &=& p(t,y,x),\nonumber\\
p(t+s,x,y) &=& \int_M p(t,x,z) \,p(s,z,y)\,dz, \nonumber\\
\int_M p(t,x,y)\,dy &\leq& 1 .
\label{eq:ptotal}
\end{eqnarray}

The heat kernel has the following stochastic interpretation.
For $x \in M$ and $U \subset M$ open, $\int_U p(t,x,y)\,dy$ is the probability that a random path emanating from $x$ lies in $U$ at time $t$.
Thus if we have strict inequality in \eqref{eq:ptotal}, then there is a positive probability that a random path will reach infinity in finite time $t$.
This motivates the following

\begin{dfn}
A geodesically complete connected Riemannian manifold is called \emph{stochastically complete} if $\int_M p(t,x,y)\,dy = 1$ for some (or equivalently all) $t>0$ and $x\in M$.
\end{dfn}

The concept of stochastic completeness can also be considered for geodesically incomplete manifolds but we will not need this.

Various sufficient geometric criteria for stochastic completeness of geodesically complete manifolds are known.
Yau \cite[Cor.~2]{Y} showed that if the Ricci curvature is bounded from below, then $M$ has no non-zero bounded eigenfunctions of $\Delta$ for eigenvalues $\lambda \gg 0$.
By \cite[Thm.~6.2, Crit.~3]{G} this shows that $M$ is stochastically complete.

Grigor'yan \cite[Thm.~9.1]{G} has a very nice criterion in terms of volume growth.
For any $x\in M$ denote the closed ball of radius $r>0$ about $x$ by $B(x,r)$.
We write $V(x,r) := \vol (B(x,r))$ and $S(x,r) := \area (\partial B(x,r))$.
Here $\vol$ denotes the $n$-dimensional volume and $\area$ the $(n-1)$-dimensional volume.
Now Grigor'yan's criterion says that if 
\begin{equation}
\int^\infty \frac{r\,dr}{\log V(x,r)} \ep=\ep \infty
\label{eq:GrigsCrit}
\end{equation}
for some $x\in M$, then $M$ is stochastically complete.
Note that this criterion can be applied if $V(x,r) \leq \exp(C\cdot r^2)$ for some $C>0$ and all $r\geq r_0$.

There is a particularly simple class of spherically symmetric manifolds for which one can study geometric properties of stochastically complete manifolds rather explicitly.
They are sometimes called ``model manifolds'' in this context and they arise as follows.
Let $f:[0,\infty) \to \R$ be a smooth function such that $f(0)=0$, $f'(0)=1$, and $f(t)>0$ for $t>0$.
Then we call $\R^n$ equipped with the metric $g=dr^2+f(r)^2g_{S^{n-1}}$ a \emph{model manifold}.
Here $r=|x|$ is the distance from the origin $o\in\R^n$ and $g_{S^{n-1}}$ is the standard metric of $S^{n-1}$.
For example, Euclidean space and hyperbolic space are model manifolds with $f(r)=r$ and $f(r)=\sinh(r)$ respectively.
It is not too hard to show \cite[Prop.~3.2]{G} that a model manifold is stochastically complete if and only if 
\begin{equation}
\int^\infty \frac{V(o,r)}{S(o,r)}\, dr \ep=\ep \infty .
\nonumber
\end{equation}

\begin{example}\label{ex:expralpha}
Let $\alpha\in \R$ and let $f(r)= r^{(\alpha-1)/(n-1)} \exp\left(\frac{r^\alpha}{n-1}\right)$ for $r\geq1$.
Then for $r\geq1$ we have $S(o,r) = C_1 \cdot f(r)^{n-1} = C_1 \cdot r^{\alpha-1} \exp\left(r^\alpha\right)$ and $V(o,r) = C_2 + \int_{1}^r S(o,\rho)\,d\rho = C_2 + C_3 \cdot \exp\left(r^\alpha\right)$.
Hence
$$
\int_1^\infty \frac{V(o,r)}{S(o,r)}\, dr 
\ep=\ep
\int_1^\infty \frac{C_2\cdot\exp(r^{-\alpha})+C_3}{C_1\cdot r^{\alpha-1}} dr
\ep=\ep
\infty
$$
if and only if $\alpha \leq 2$.
This shows that Grigor'yan's criterion \eqref{eq:GrigsCrit} is quite sharp.
\end{example}

It should be noted that the much stronger condition
$$
\int^\infty \frac{dr}{S(o,r)} \ep=\ep \infty 
$$
is equivalent to Brownian motion on the model manifold being recurrent.
Lyons and Sullivan \cite[Sec.~6]{LS} and Grigor'yan \cite{G1,G2} independently showed that for a general geodesically complete manifold $M$ the condition
$$
\int^\infty \frac{dr}{S(x,r)} \ep=\ep \infty 
$$
for some $x\in M$ implies recurrence of Brownian motion on $M$.
However, on non-model manifolds this condition is not necessary for recurrence of the Brownian motion as can be shown by examples \cite[Example~7.3]{G}.
Grigor'yan asked \cite[Problem~9]{G} if similarly on a general geodesically complete manifold $M$ the condition
\begin{equation}
\int^\infty \frac{V(x,r)}{S(x,r)}\, dr \ep=\ep \infty .
\label{eq:ConjCrit}
\end{equation}
for some $x\in M$ is sufficient for stochastic completeness.
Sometimes this is formulated as a conjecture \cite[Remark on p.~40]{PRS}.
The main result of the present paper is the construction of counter-examples to this conjecture.

\begin{thm}\label{thm:main}
In any dimension $n\geq 2$ there exists a geodesically complete but stochastically incomplete  connected Riemannian manifold $M$ such that for some $x\in M$ the volume growth condition \eqref{eq:ConjCrit} holds.
\end{thm}

Thus the analog to the result of Lyons, Sullivan, and Grigor'yan for stochastic completeness does not hold.

\section{The weak Omori-Yau maximum principle}

As a useful tool we recall the \emph{weak Omori-Yau maximum principle}.
It says that for each $u \in C^2(M)$ with $u^* := \sup_Mu < \infty$ there exists a sequence $x_k\in M$ such that
\begin{equation}
\lim_{k\to \infty} u(x_k) \ep=\ep u^* ,
\label{OY1}
\end{equation}
\begin{equation}
\limsup_{k\to\infty} \Delta u(x_k) \ep\leq\ep 0 .
\label{OY2}
\end{equation}
It is a theorem by Pigola, Rigoli and Setti \cite{PRS1},\cite[Thm.~3.1]{PRS} that the validity of the weak Omori-Yau maximum principle is equivalent to $M$ being stochastically complete.
In other words, a stochastically incomplete manifold is characterized by the existence of a function $u \in C^2(M)$ with $u^* := \sup_Mu < \infty$ such that for any sequence $x_k\in M$ satisfying \eqref{OY1} we have
\begin{equation}
\limsup_{k\to\infty} \Delta u(x_k) \ep > \ep 0 .
\label{OY3}
\end{equation}
We will call such a function \viol.
It is clear from the definition that if $u$ is \viol\ and $v\in C^2(M)$ coincides with $u$ outside a compact subset $K \subset M$ and $v<u^*$ on $K$, then $v$ is \viol\ as well.

\begin{example}\label{ex:OY}
Let $f(r)= r^{(\alpha-1)/(n-1)} \exp\left(\frac{r^\alpha}{n-1}\right)$ for $r\geq1$ be as in Example~\ref{ex:expralpha} with $\alpha>2$.
We know that the corresponding model manifold is stochastically incomplete.
To exhibit a \viol\ function choose $\beta>0$ such that $\alpha - \beta > 2$.
Now let $u$ be a smooth function on the model manifold depending on $r$ only such that $u(r) = 1-r^{-\beta}$ for $r\geq R_1$ and $u<1$ everywhere.
Then $u^*=1$.
On a model manifold the Laplace operator takes the form 
$$
\Delta \ep = \ep
\frac{\partial^2}{\partial r^2} + (n-1)\frac{f'(r)}{f(r)}\frac{\partial}{\partial r} + \frac{1}{f(r)^2}\Delta_S
$$
where $\Delta_S$ is the Laplace-Beltrami operator on the standard sphere $S^{n-1}$.
Hence for $r\geq R_1$
$$
\Delta u
\ep = \ep
u''(r) + (n-1)\frac{f'(r)}{f(r)} u'(r)
\ep = \ep
\beta\left( \alpha\, r^{\alpha-\beta-2} + (\alpha-\beta-2)\, r^{-\beta-2}\right).
$$
This goes to $\infty$ as $r\to\infty$.
Since for any sequence $r_k$ such that $u(r_k) \to 1$ we must have $r_k \to \infty$ we see that $u$ is \viol.
\end{example}

\section{Stochastic completeness and connected sums}

Our examples will be constructed as connected sums.
Hence we have first to examine to what extent stochastic completeness is preserved under this operation.

\begin{lemma}\label{lem:ZusSumme}
Let $M_1$ and $M_2$ be geodesically complete Riemannian manifolds of equal dimension.
Let $K \subset M_1 \sharp M_2$, $K_1 \subset M_1$, and $K_2 \subset M_2$ be compact subsets such that $M_1 \sharp M_2 \setminus K$ is isometric to the disjoint union of $M_1 \setminus K_1$ and $M_2 \setminus K_2$.

Then $M_1 \sharp M_2$ is stochastically complete if and only if  $M_1$ and $M_2$ are stochastically complete.
\end{lemma}

\begin{center}
{
\begin{pspicture}(-1,-2.3)(12.62,2.6)
\psset{unit=8.5mm}

\pscustom[fillstyle=solid,fillcolor=gray]{
\psbezier[linewidth=0.04](0.0,2.91)(0.0,2.11)(5.380851,0.8497826)(5.36,-0.15)(5.3391495,-1.1497827)(1.02,-1.53)(0.02,-2.91)
}
\pscustom[linewidth=0pt,linecolor=white,fillstyle=solid,fillcolor=white]{
\psline(-1,3)(3.87,3)(3.87,-3)(-1,-3)
}
\pscustom[linewidth=0pt,linecolor=gray,fillstyle=solid,fillcolor=gray]{
\psellipticarc[linewidth=0.04,dimen=outer](3.87,-0.08)(0.55,1.03){90}{270}
}
\psellipse[linewidth=0pt,dimen=outer,linecolor=gray,fillstyle=solid,fillcolor=gray](3.87,-0.08)(0.55,1.03)

\psbezier[linewidth=0.04](0.0,2.91)(0.0,2.11)(5.380851,0.8497826)(5.36,-0.15)(5.3391495,-1.1497827)(1.02,-1.53)(0.02,-2.91)
\psellipticarc[linewidth=0.04,dimen=outer](3.87,-0.08)(0.55,1.03){90}{270}
\psellipticarc[linewidth=0.04,dimen=outer,linestyle=dashed](3.87,-0.08)(0.55,1.03){270}{90}
\rput(2,-0.1){$M_1$}
\rput(4.4,-0.1){\psframebox*[framearc=.7]{$K_1$}}

\pscustom[fillstyle=solid,fillcolor=gray]{
\psbezier[linewidth=0.04](12.6,1.89)(11.8,1.89)(6.2312984,0.84958804)(6.26,-0.15)(6.2887015,-1.149588)(11.76,-1.85)(12.58,-1.79)
}
\pscustom[linewidth=0pt,linecolor=white,fillstyle=solid,fillcolor=white]{
\psline(14,2)(7.59,2)(7.59,-2)(14,-2)
}
\pscustom[linewidth=0pt,linecolor=white,fillstyle=solid,fillcolor=white]{
\psellipticarc[linewidth=0.04,dimen=outer](7.59,-0.11)(0.47,0.86){90}{270}
}
\psellipse[linewidth=0pt,dimen=outer,linecolor=lightgray,fillstyle=solid,fillcolor=lightgray](7.59,-0.11)(0.47,0.86)
\psbezier[linewidth=0.04](9.84,0.27)(9.84,-0.53)(11.62,-0.55)(11.62,0.25)
\psbezier[linewidth=0.04](11.420144,-0.18283969)(11.407777,0.36940283)(10.08749,0.36746314)(10.099856,-0.18477936)
\psbezier[linewidth=0.04](12.6,1.89)(11.8,1.89)(6.2312984,0.84958804)(6.26,-0.15)(6.2887015,-1.149588)(11.76,-1.85)(12.58,-1.79)
\psellipticarc[linewidth=0.04,dimen=outer](7.59,-0.11)(0.47,0.86){90}{270}
\psellipticarc[linewidth=0.04,dimen=outer,linestyle=dashed](7.59,-0.11)(0.47,0.86){270}{90}
\rput(9,-0.1){$M_2$}
\rput(6.74,-0.1){\psframebox*[framearc=.7]{$K_2$}}
\end{pspicture} 
}
\end{center}

\begin{center}
{
\begin{pspicture}(-1,-2.93)(12.62,2.3)
\psset{unit=8.5mm}
\pscustom[fillstyle=solid,fillcolor=gray]{
\psbezier[linewidth=0.04](0.0,2.91)(0.0,2.11)(5.380851,0.8497826)(5.36,-0.15)(5.3391495,-1.1497827)(1.02,-1.53)(0.02,-2.91)
}
\pscustom[linewidth=0pt,linecolor=white,fillstyle=solid,fillcolor=white]{
\psline(-1,3)(3.87,3)(3.87,-3)(-1,-3)
}
\pscustom[linewidth=0pt,linecolor=gray,fillstyle=solid,fillcolor=gray]{
\psellipticarc[linewidth=0.04,dimen=outer](3.87,-0.08)(0.55,1.03){90}{270}
}
\psellipse[linewidth=0pt,dimen=outer,linecolor=gray,fillstyle=solid,fillcolor=gray](3.87,-0.08)(0.55,1.03)

\psbezier[linewidth=0.04](0.0,2.91)(0.0,2.11)(5.380851,0.8497826)(5.36,-0.15)(5.3391495,-1.1497827)(1.02,-1.53)(0.02,-2.91)
\psellipticarc[linewidth=0.04,dimen=outer](3.87,-0.08)(0.55,1.03){90}{270}
\psellipticarc[linewidth=0.04,dimen=outer,linestyle=dashed](3.87,-0.08)(0.55,1.03){270}{90}

\pscustom[fillstyle=solid,fillcolor=gray]{
\psbezier[linewidth=0.04](12.6,1.89)(11.8,1.89)(6.2312984,0.84958804)(6.26,-0.15)(6.2887015,-1.149588)(11.76,-1.85)(12.58,-1.79)
}
\pscustom[linewidth=0pt,linecolor=white,fillstyle=solid,fillcolor=white]{
\psline(14,2)(7.59,2)(7.59,-2)(14,-2)
}
\pscustom[linewidth=0pt,linecolor=white,fillstyle=solid,fillcolor=white]{
\psellipticarc[linewidth=0.04,dimen=outer](7.59,-0.11)(0.47,0.86){90}{270}
}
\psellipse[linewidth=0pt,dimen=outer,linecolor=lightgray,fillstyle=solid,fillcolor=lightgray](7.59,-0.11)(0.47,0.86)

\psbezier[linewidth=0.04](9.84,0.27)(9.84,-0.53)(11.62,-0.55)(11.62,0.25)
\psbezier[linewidth=0.04](11.420144,-0.18283969)(11.407777,0.36940283)(10.08749,0.36746314)(10.099856,-0.18477936)
\psbezier[linewidth=0.04](12.6,1.89)(11.8,1.89)(6.2312984,0.84958804)(6.26,-0.15)(6.2887015,-1.149588)(11.76,-1.85)(12.58,-1.79)
\psellipticarc[linewidth=0.04,dimen=outer](7.59,-0.11)(0.47,0.86){90}{270}
\psellipticarc[linewidth=0.04,dimen=outer,linestyle=dashed](7.59,-0.11)(0.47,0.86){270}{90}
\pscustom[linewidth=0.04,linecolor=gray,fillstyle=solid,fillcolor=gray]{
\psbezier[linewidth=0.04](4.6,0.31)(5.34,-0.23)(6.46,-0.15)(7.04,0.23)
\psbezier[liftpen=1,linewidth=0.04](7.04,-0.49)(6.6,-0.23)(5.36,-0.07)(4.68,-0.49)
}
\psbezier[linewidth=0.04](4.6,0.31)(5.34,-0.23)(6.46,-0.15)(7.04,0.23)
\psbezier[liftpen=1,linewidth=0.04](7.04,-0.49)(6.6,-0.23)(5.36,-0.07)(4.68,-0.49)
\rput(2,-0.1){$M_1 \sharp M_2$}
\rput(4.4,-0.1){\psframebox*[framearc=.7]{$K$}}
\end{pspicture} 
}

\emph{Fig.~1}
\end{center}

\begin{proof}
Suppose that $M_1$ is stochastically incomplete.
Then there exists a \viol\ function $u \in C^2(M_1)$.
By adding a constant if necessary we can w.l.o.g.\ assume that $u^*>0$.
Let $\chi \in C^\infty(M_1)$ be a function satisfying $0 \leq \chi \leq 1$ on $M_1$, $\chi \equiv 0$ on $K_1$ and $\chi\equiv 1$ outside a compact set.
Put $v:=\chi\cdot u$.
Then $v\in C^2(M_1)$ coincides with $u$ outside a compact set and $v<u^*$ on this compact set.
Thus $v$ is \viol\ as well.
Since $v$ vanishes on $K_1$ we can extend it by zero and regard it as a function on $M_1 \sharp M_2$.
Thus we have a \viol\ function on the connected sum which shows that $M_1 \sharp M_2$ is stochastically incomplete.

Conversely, let $M_1 \sharp M_2$ be stochastically incomplete.
Let $u\in C^2(M_1 \sharp M_2)$ be a \viol\ function and assume again that $0<u^* < \infty$.
Let $\chi_1\in C^\infty(M_1 \sharp M_2)$ be a function satisfying $0 \leq \chi_1 \leq 1$ on all of $M_1 \sharp M_2$, $\chi_1 \equiv 0$ on $K \cup M_2$ and $\chi_1 \equiv 1$ outside a compact subset of $M_1$.
Define $\chi_2$ similarly by interchanging the roles of $M_1$ and $M_2$.
Put $u_j := \chi_j \cdot u$.
Then $u=u_1 +u_2$ outside a compact subset of $M_1 \sharp M_2$ and $u_1+u_2 < u^*$ everywhere.
By similar reasoning as above $u=u_1 +u_2$ is \viol, hence $u_1$ or $u_2$ is \viol\ as well.
Since $u_j$ can be considered as a function on $M_j$ we conclude that $M_1$ or $M_2$ must be stochastically incomplete.
\end{proof}

\begin{rem}
Another criterion for stochastic incompleteness which can be used for an easy proof of Lemma~\ref{lem:ZusSumme} is that $M$ is $\lambda$-massive \cite[Thm.~6.2]{G}.
By \cite[Prop.~6.1]{G} $\lambda$-massiveness of a subset of a Riemannian manifold is preserved by enlarging the subset and also by subtracting a compact subset.
Hence if $M_1$ is stochastically incomplete, then $\Omega_1 = M_1 \setminus K_1$ is $\lambda$-massive.
Thus $M = M_1 \sharp M_2=\Omega_1 \cup \Omega_2 \cup K \supset \Omega_1$ is $\lambda$-massive and therefore stochastically incomplete.
The converse implication is proved similarly.
\end{rem}

\section{Construction of the counter-example}

To construct the counter-examples and prove Theorem~\ref{thm:main} we pick a geodesically complete but stochastically incomplete Riemannian manifold $M_1$.
Specifically, we may take a model manifold as in Example~\ref{ex:expralpha} with $\alpha >2$.
To prepare for the connected sum we fix a compact subset $K_1 \subset M_1$ with non-empty interior and remove a small open ball from the interior of $K_1$.
We obtain a manifold $\widehat M_1$ with boundary diffeomorphic to $S^{n-1}$.
After a deformation of the Riemannian metric inside $K_1$ we can assume that near the boundary the metric is of product form $dr^2 + C_1^2 \cdot g_{S^{n-1}}$ where the scaling factor $C_1>0$ is chosen such that the (intrinsic) diameter of the boundary is $1/8$.
Fix $q_1 \in \partial \widehat M_1$ and put 
$$
S_1(r) \ep:=\ep \area (\partial B^{\widehat M_1}(q_1,r))
$$
and 
$$
F(r) \ep:=\ep \max_{\rho\in[0,r]}S_1(\rho).
$$
Then $F$ is a monotonically increasing function.
Next we choose a smooth function $V:[0,\infty) \to \R$ such that
\begin{itemize}
\item 
$V(0)\ep=\ep 0$
\item
$S(r) \ep:=\ep V'(r)\ep>\ep 0$ for all $r\in [0,\infty)$
\item
$V(k) \ep\geq\ep F(k+1)$ for all $k=1,2,3,\ldots$
\item
$S$ is constant on all intervals $[k+\frac18,k+\frac34]$, $k=1,2,3,\ldots$
\end{itemize}
The model manifold with warping function $f(r) := \sqrt[n-1]{S(r)/\omega_{n-1}}$ has $V(o,r) = V(r)$ and $S(o,r) = S(r)$.
Deform $f$ near $0$ such that $f(r) = C_1$ for $r$ near $0$ and $\int_0^1 f(r)^{n-1}dr$ remains unchanged.

\begin{center}
{
\begin{pspicture}(0,-3.3)(5,3)
\psset{unit=4mm}

\psbezier[linewidth=0.04](0.06,-7.48)(0.82,-6.68)(2.32,-3.08)(3.16,-3.08)(4.0,-3.08)(3.66,-3.08)(4.66,-3.08)(5.66,-3.08)(5.06,-3.08)(5.86,-3.08)(6.66,-3.08)(5.92,2.5)(6.68,2.52)(7.44,2.54)(7.140114,2.5049007)(8.14,2.52)(9.139886,2.5350993)(8.24,2.5)(9.32,2.52)(10.4,2.54)(9.76,6.5)(10.02,7.0)
\psline[linewidth=0.04cm]{->}(0.06,-7.5)(11.48,-7.5)
\psline[linewidth=0.04cm]{->}(0.06,-7.5)(0.06,7.5)
\psline[linewidth=0.04cm](3.26,-7.5)(3.26,-7.3)
\psline[linewidth=0.04cm](6.46,-7.5)(6.46,-7.3)
\psline[linewidth=0.04cm](9.66,-7.5)(9.66,-7.3)
\psbezier[linewidth=0.04,linestyle=dashed](2.12,-4.1)(1.6,-4.64)(2.18,-6.3)(1.28,-6.3)(0.38,-6.3)(1.1209902,-6.3)(0.04,-6.3)
\rput(-0.5,-6.3){$C_1$}
\rput(4.6,-2.55){$f$}
\rput(3.26,-6.9){$1$}
\rput(6.46,-6.9){$2$}
\rput(9.66,-6.9){$3$}
\end{pspicture} 
}

\emph{Fig.~2}
\end{center}

Let $\widehat M_2$ be the manifold $[0,\infty) \times S^{n-1}$ with the Riemannian metric $dr^2 + f(r)^2\cdot g_{S^{n-1}}$.
Then $\widehat M_2$ is a manifold with boundary diffeomeomorphic to $S^{n-1}$ such that the diameter of $\partial \widehat M_2$ is $1/8$.
Furthermore, for all $r\geq1$,
$$
V(r) \ep=\ep \vol(\{x\in \widehat M_2\,|\, d(x,\partial \widehat M_2) \leq r\})
$$
and
$$
S(r) 
\ep=\ep 
\area(\{x\in \widehat M_2\,|\, d(x,\partial \widehat M_2) = r\})
\ep =\ep 
\area(\partial\{x\in \widehat M_2\,|\, d(x,\partial \widehat M_2) \leq r\}).
$$
Pick $q_2 \in \partial \widehat M_2$.
Put $V_2(r) := \vol(B^{\widehat M_2}(q_2,r))$ and $S_2(r) := \area (\partial B^{\widehat M_2}(q_2,r))$.
By the triangle inequality we have for all $r\geq1$
$$
\{x\in \widehat M_2\,|\, d(x,\partial \widehat M_2) \leq r-1/8\} 
\ep\subset\ep 
B^{\widehat M_2}(q_2,r) 
\ep\subset\ep 
\{x\in \widehat M_2\,|\, d(x,\partial \widehat M_2) \leq r\}
$$
and hence
$$
V(r-1/8) \ep\leq\ep V_2(r) \ep\leq\ep V(r) .
$$
Now we glue $\widehat M_1$ and $\widehat M_2$ along the boundary such that $q_1$ and $q_2$ get identified to one point $q$.
This yields a smooth and geodesically complete Riemannian metric on $M = M_1 \sharp M_2$.
Since $M_1$ is stochastically incomplete, so is $M$ by Lemma~\ref{lem:ZusSumme}.
It remains to show that
\begin{equation}
\int_1^\infty \frac{V(q,r)}{S(q,r)}\, dr 
\ep=\ep 
\int_1^\infty \frac{V_1(r)+V_2(r)}{S_1(r)+S_2(r)}\,dr 
\ep=\ep 
\infty .
\label{eq:ConjCrit2}
\end{equation}
For this purpose we estimate $\frac{S_1(r)+S_2(r)}{V_1(r)+V_2(r)}$ for $r\in[k+\frac12,k+\frac34]$, $k\in\N$.
Namely,
\begin{eqnarray*}
\frac{S_1(r)+S_2(r)}{V_1(r)+V_2(r)}
&\leq&
\frac{S_1(r)+S_2(r)}{V_2(r)} \\
&\leq&
\frac{F(k+1)+S_2(r)}{V_2(r)} \\
&\leq&
\frac{F(k+1)}{V(r-\frac18)} + \frac{S_2(r)}{V_2(r)} \\
&\leq&
\frac{F(k+1)}{V(k)} + \frac{S_2(r)}{V_2(r)} \\
&\leq&
1 + \frac{S_2(r)}{V_2(r)} .
\end{eqnarray*}
By the Cauchy-Schwarz inequality we find
$$
\frac{1}{16} 
\ep=\ep 
\left(\int_{k+\frac12}^{k+\frac34}1\,dr\right)^2
\ep\leq\ep
\left(\int_{k+\frac12}^{k+\frac34}\frac{S_1(r)+S_2(r)}{V_1(r)+V_2(r)}\,dr\right) \cdot
\left(\int_{k+\frac12}^{k+\frac34}\frac{V_1(r)+V_2(r)}{S_1(r)+S_2(r)}\,dr\right) ,
$$
hence
\begin{eqnarray}
16\int_{k+\frac12}^{k+\frac34}\frac{V_1(r)+V_2(r)}{S_1(r)+S_2(r)}\,dr
&\geq&
\left(\int_{k+\frac12}^{k+\frac34}\frac{S_1(r)+S_2(r)}{V_1(r)+V_2(r)}\,dr\right)^{-1} \nonumber\\
&\geq&
\left(\int_{k+\frac12}^{k+\frac34}\left(1 + \frac{S_2(r)}{V_2(r)}\right)\,dr\right)^{-1} \nonumber\\
&=&
\left(\frac14 + \int_{k+\frac12}^{k+\frac34}\left(\frac{d}{dr}\log(V_2(r))\right)\,dr\right)^{-1} \nonumber\\
&=&
\left(\frac14 + \log(V_2(k+3/4)) - \log(V_2(k+1/2))\right)^{-1} \nonumber\\
&\geq&
\left(\frac14 + \log(V(k+3/4)) - \log(V(k+1/2-1/8))\right)^{-1} \nonumber\\
&=&
\left(\frac14 + \int_{k+\frac38}^{k+\frac34} \frac{S(r)}{V(r)}\, dr\right)^{-1} .
\label{eq:est1}
\end{eqnarray}
Since $S=V'$ is constant on $[k+\frac18,k+\frac34]$ we have for $r\in[k+\frac38,k+\frac34]$ that $S(r) = S(k+\frac18)$ and $V(r) \geq S(k+\frac18)\cdot (3/8-1/8)=S(k+\frac18)/4$.
Thus 
$$
\int_{k+\frac38}^{k+\frac34} \frac{S(r)}{V(r)}\, dr 
\ep\leq\ep 
\frac14 \cdot \left(\frac34-\frac38\right) 
\ep=\ep 
\frac{3}{32} .
$$
Plugging this into \eqref{eq:est1} yields
$$
\int_{k+\frac12}^{k+\frac34}\frac{V_1(r)+V_2(r)}{S_1(r)+S_2(r)}\,dr
\ep\geq\ep \frac{2}{11} .
$$
Summation over $k$ gives
$$
\int_1^\infty \frac{V_1(r)+V_2(r)}{S_1(r)+S_2(r)}\,dr \ep=\ep \infty
$$
as desired.
This concludes the construction of the counter-example and the proof of Theorem~\ref{thm:main}.

\section{Concluding remarks}

\begin{rem}
The examples constructed in the previous section have (at least) two ends.
One may ask whether or not one can find examples with only one topological end.\footnote{We thank B.~Wilking for bringing up this question.}
Indeed, this is possible.

One starts with an example $M=M_1 \sharp M_2$ with two ends as constructed above.
Let $u\in C^2(M)$ be a \viol\ function vanishing on the second end $\widehat M_2$ and such that $0<u^*<\infty$ as constructed in the proof of Lemma~\ref{lem:ZusSumme}.
Choose a sequence of points $x_k\in \widehat M_1$ in the first end of $M$ satisfying \eqref{OY1} and \eqref{OY3}.
Then $r_k := d(q,x_k) \to \infty$ as $k\to \infty$.
Now pick a monotonically increasing sequence of numbers $R_j>0$ such that $R_j\to \infty$ as $j\to\infty$ and $r_k \neq R_j$ for all $k$ and $j$.
We choose $\eps_j>0$ so small that the intervals $(R_j-\eps_j,R_j+\eps_j)$ are pairwise disjoint, such that $r_k \not\in (R_j-\eps_j,R_j+\eps_j)$ for all $k$ and $j$ and such that 
\begin{equation}
\sum_{j=1}^\infty \int_{R_j-\eps_j}^{R_j+\eps_j} \frac{V(q,r)}{S(q,r)}\, dr 
\ep<\ep \infty .
\label{eq:stoerendl}
\end{equation}
The minimal geodesics from $q$ to $\{x\in \widehat M_1\,|\, d(q,x)=R_j+\eps_j\}$ do not cover all of $B(q,R_j+\eps_j) \cap \widehat M_1$.
The complement is a non-empty open ``wedge'' whose boundary intersects $\{x\in \widehat M_1\,|\, d(q,x)=R_j+\eps_j\}$ at a point opposite to $q$ on $S^{n-1}$.
\begin{center}
\begin{pspicture}(-7,-3)(7,3)
\psset{unit=10mm}

\psellipticarc(-3,0)(1,2){80}{280}
\psellipticarc[linewidth=0.3pt,linestyle=dashed](-3,0)(1,2){280}{80}

\psecurve(4,-3)(3.5,0)(1,1.03)(-4,0)(-5,-1)
\psecurve(4,3)(3.5,0)(1,-1.03)(-4,0)(-5,1)
\psline[linecolor=white,fillstyle=solid,fillcolor=white](1,2)(2.64,2)(2.64,-2)(1,-2)
\psecurve[linewidth=0.3pt,linestyle=dashed](4,-3)(3.5,0)(1,1.03)(-4,0)(-5,-1)
\psecurve[linewidth=0.3pt,linestyle=dashed](4,3)(3.5,0)(1,-1.03)(-4,0)(-5,1)

\psellipse(3,0)(0.5,1)
\psdot(3.5,0)
\rput(3.8,0){$q$}

\psecurve(4,1.3)(3,1)(-4,2.64)(-5,4)
\psecurve(4,-1.3)(3,-1)(-4,-2.64)(-5,-4)

\psdot(-4,0)

\psarc(-3.3,0){0.2}{45}{7}
\psdot(-3.3,0)

\rput(-2.6,0.2){$y_j$}
\psline[linewidth=0.3pt]{->}(-2.8,0.2)(-3.25,0.02)
\rput(-5.5,2){$\{d(x,q)=R_j+\eps_j\}$}
\psline[linewidth=0.3pt]{->}(-4,2)(-3.5,1.8)
\rput(-4,-2){$\widehat M_1$}
\end{pspicture}
\emph{Fig.~3}
\end{center}
Choose points $y_j \in \widehat M_1$ with $d(q,y_j) = R_j$ and $\delta_j\in (0,\eps_j/2)$ so small that $B(y_j,\delta_j)$ is contained in this wedge.
Moreover choose $z_j \in \widehat M_2$ with $d(q,z_j)= R_j$.

We remove the balls $B(y_j,\delta_j)$ and $B(z_j,\delta_j)$ from $M$ and glue in handles $H_j$ diffeomorphic to $S^{n-1}\times [0,1]$.
We denote the resulting manifold by $\widetilde M$.
The handles connect the two ends of $M$ outside each compact set so that $\widetilde M$ has only one topological end.

We choose the metric on the handles $H_j$ such that $\vol(H_j) = \vol (B(y_j,\delta_j)) + \vol (B(z_j,\delta_j))$, such that minimal geodesics through $H_j$ joining two points on $\partial B(y_j,\delta_j)$ (or two points on $\partial B(z_j,\delta_j)$) are no shorter than those through $B(y_j,\delta_j)$ (or $B(z_j,\delta_j)$ resp.) and such that we obtain a smooth metric on $\widetilde M$.
To see that such metrics exist on $H_j$ we first look at the case that $B(y_j,\delta_j)$ and $B(z_j,\delta_j)$ are isometric to Euclidean balls.
Then the metric can be chosen such that $H_j$ is a cylinder flattened near the two boundary components.
The flattening ensures that the metric extends smoothly to $\widetilde M$, the height of the cylinder can be chosen such that the volume is right and the condition on the length of geodesics is also fulfilled.

\begin{center}
\begin{pspicture}(-6,-3)(6,3.5)
\psset{unit=10mm}

\psellipse(-3,-1.5)(2,1)
\psellipse(-3,2)(2,1)
\psdot(-3,-1.5)
\rput(-2.1,-1.5){$B(y_j,\delta_j)$}
\psdot(-3,2)
\rput(-2.1,2){$B(z_j,\delta_j)$}
\psline[linewidth=0.5pt](-4.7,-1)(-3,-2.5)
\psdots(-4.7,-1)(-3,-2.49)

\psecurve(1.2,2)(1.5,2)(1.8,1.8)(2.05,0.25)(1.8,-1.3)(1.5,-1.5)(1.2,-1.5)
\psecurve(4.8,2)(4.5,2)(4.2,1.8)(3.95,0.25)(4.2,-1.3)(4.5,-1.5)(4.8,-1.5)
\pspolygon[fillstyle=solid,fillcolor=white,linecolor=white](1.5,1.68)(2.1,1.68)(2.1,1.15)(1.5,1.15)
\pspolygon[fillstyle=solid,fillcolor=white,linecolor=white](4.5,1.68)(3.9,1.68)(3.9,1.15)(4.5,1.15)
\psellipse(3,2)(2,1)
\psellipse(3,2)(1.3,0.65)
\psellipse[linewidth=0.3pt,linestyle=dashed](3,-1.5)(2,1)
\psellipticarc(3,-1.5)(2,1){120}{60}
\psellipse[linewidth=0.3pt,linestyle=dashed](3,0.25)(0.95,0.47)
\psellipticarc(3,0.25)(0.95,0.47){180}{0}
\psdots(1.3,-1)(3,-2.49)
\psecurve[linewidth=0.5pt](-0.4,0.5)(1.3,-1)(1.81,-1.25)(3,-0.7)
\psecurve[linewidth=0.5pt](1.3,-1.8)(1.9,-1.1)(3,-2.49)(4.7,-3.99)
\rput(4.3,0.25){$H_j$}

\end{pspicture}
\emph{Fig.~4}
\end{center}

This construction is robust under slight perturbations of the metrics.
Hence, in the general case of curved balls $B(y_j,\delta_j)$ and $B(z_j,\delta_j)$ we choose $\delta_j$ so small that the balls are sufficiently close to Euclidean balls so that the same construction still works.

With these choices we have
$$
\widetilde V(q,r) = V(q,r)
\mbox{ and }
\widetilde S(q,r) = S(q,r)
$$
for all $r>0$ not lying in any of the intervals $[R_j-\eps_j,R_j+\eps_j]$.
Here $\widetilde V$ and $\widetilde S$ denote the volumes of the balls and of their boundaries in $\widetilde M$.
Therefore, by \eqref{eq:stoerendl},
\begin{eqnarray*}
\int_0^\infty \frac{\widetilde V(q,r)}{\widetilde S(q,r)}\, dr
&\geq&
\int_{(0,\infty) \setminus \cup_{j=1}^\infty [R_j-\eps_j,R_j+\eps_j]} \frac{\widetilde V(q,r)}{\widetilde S(q,r)}\, dr \\
&=&
\int_{(0,\infty) \setminus \cup_{j=1}^\infty [R_j-\eps_j,R_j+\eps_j]} \frac{V(q,r)}{S(q,r)}\, dr \\
&=&
\int_0^\infty \frac{V(q,r)}{S(q,r)}\, dr - \sum_{j=1}^\infty \int_{R_j-\eps_j}^{R_j+\eps_j} \frac{V(q,r)}{S(q,r)}\, dr \\
&=&
\infty .
\end{eqnarray*}
In order to see that $\widetilde M$ is stochastically incomplete, we construct a \viol\ function $\widetilde v \in C^2(\widetilde M)$.
We choose a cut-off function $\chi\in C^\infty(M)$ with $0\leq \chi \leq 1$ everywhere, $\chi \equiv 1$ outside the pairwise disjoint balls $B(y_j,\eps_j/2)$, and $\chi_j \equiv 0$ on the smaller balls $B(y_j,\delta_j)$.
Put $v:= \chi \cdot u \in C^2(M)$.
Since $v\leq u^*$ everywhere and $v=u$ on neighborhoods of the $x_k$ we see that $v$ is \viol.
We restrict $v$ to $M$ minus the $\delta_j$-balls and extend it by zero over the handles.
This yields a \viol\ function $\widetilde v$ on $\widetilde M$.
\end{rem}

\begin{rem}
Conversely, one may also ask if on a general geodesically complete manifold $M$ the condition
\begin{equation}\label{eq:VSendl}
\int^\infty \frac{V(x,r)}{S(x,r)}\, dr \ep<\ep \infty
\end{equation}
for some $x\in M$ implies stochastic incompleteness.
But this is false too as we will demonstrate by a counter-example.
We start the construction with a modification of \cite[Ex.~7.3]{G}.
Choose positive smooth functions $S_1, S_2 : (0,\infty) \to \R$ with the following properties:
\begin{itemize}
\item[(P1)]
$S_1(r) = S_2(r) = 2\pi r$ for $0<r\leq 1$
\item[(P2)]
$S_1(r) + S_2(r) = 3r^2\exp(r^3)$ for $r\geq 2$
\item[(P3)]
$S_1(r) = 1$ for $r\in [4k,4k+1]$, $k\in \N$
\item[(P4)]
$S_2(r) = 1$ for $r\in [4k+2,4k+3]$, $k\in \N$
\end{itemize}
Let $M_1$ and $M_2$ be the corresponding 2-dimensional model manifolds with warping functions $f_j(r) = S_j(r)/2\pi$.
Then $S_1(r) = S(o,r)$ in $M_1$ and similarly for $M_2$.
Properties (P3) and (P4) imply
$$
\int^\infty \frac{dr}{S_j(r)} \ep = \ep \infty ,
$$
hence Brownian motion is recurrent.
In particular, $M_1$ and $M_2$ are stochastically complete and, by Lemma~\ref{lem:ZusSumme}, so is the connected sum $M_1 \sharp M_2$.

Now let $V_j(r) := V(o,r)$ in $M_j$, in other words, $V_j'=S_j$ and $V_j(0)=0$.
From (P2) we conclude $V_1(r) + V_2(r) = \exp(r^3) + C$ for $r \geq 2$.
Thus 
$$
\int^\infty \frac{V_1(r)+V_2(r)}{S_1(r)+S_2(r)}\, dr \ep<\ep \infty .
$$
To construct the metric on the connected sum $M_1 \sharp M_2$ we observe that by Property (P1) the unit disk about $o$ in $M_j$ is isometric to the unit disk in Euclidean $\R^2$.
We choose a point $p_j$ at distance $\frac12$ from $o$ and remove the interior of the disk $B(p_j,1/10)$ from $M_j$.
We obtain a manifold $\widehat M_j$ with boundary diffeomorphic to $S^1$.
We change the metric on $B(p_j,2/10)$ such that it becomes a product metric near the boundary, the volume of the unit disk $B(o,1)$ after removal of the small disk and modification of the metric is the same as before, and that distances from $o$ to points in $B(o,1) \setminus B(p_j,2/10)$ are not smaller after modification than they are before.
\begin{center}
\begin{pspicture}(-7.5,-2.7)(10,2.5)
\psset{unit=10mm}
\pscircle(-3,0){2}
\psdot(-3,0)
\rput(-3,0.3){$o$}
\pscircle[fillstyle=solid,fillcolor=lightgray](-3,-1){0.4}
\psarc[fillstyle=solid,fillcolor=white](-3,-1){0.2}{47}{7}
\psline[linewidth=0.2pt,linestyle=dashed](-5,0)(-1,0)
\pswedge[linewidth=0.2pt,linestyle=dashed](-3,0){2}{245}{295}
\psecurve[linewidth=0.7pt](-4.3,-0.5)(-3.85,-1.8)(-3,-1.6)(-2.5,-0.8)
\psecurve[linewidth=0.7pt](-3.5,-0.8)(-3,-1.6)(-2.15,-1.8)(-1.7,-0.5)
\rput(-6.5,-1){$\partial B^{\widetilde M_j}(o,1)$}
\psline[linewidth=0.2pt]{->}(-5.6,-1.05)(-3.3,-1.7)
\rput(-6,-1.7){$\partial B^{M_j}(o,1)$}
\psline[linewidth=0.2pt]{->}(-5.1,-1.75)(-3.4,-2)
\psdot(-3,-1)
\rput(-2.0,-0.7){$p_j$}
\psline[linewidth=0.2pt]{->}(-2.25,-0.7)(-2.91,-0.95)

\psellipse(3,0)(2,1.5)
\psdot(3,0)
\rput(3,0.3){$o$}
\psline[linecolor=white,fillcolor=lightgray,fillstyle=solid](2.8,-0.35)(3.2,-0.35)(3.2,-1)(2.8,-1)
\psellipticarc[fillcolor=lightgray,fillstyle=solid](3,-0.75)(0.4,0.3){130}{50}
\psellipse[fillcolor=white,fillstyle=solid](3,-0.35)(0.22,0.15)
\psecurve(2.8,0)(2.8,-0.35)(2.7,-0.75)(2,-0.5)
\psecurve(3.2,0)(3.2,-0.35)(3.3,-0.75)(4,-0.5)
\end{pspicture}
\emph{Fig.~5}

\end{center}
Since the disk $B(p_j,2/10)$ on which all modifications were performed is entirely contained in a half-plane with boundary containing $o$, the distance spheres from $o$ in $M_j$ and  in $\widetilde M_j$ coincide on at least one hemi-sphere.
Where they differ $\partial B^{\widetilde M_j}(o,r)$ lies inside $\partial B^{M_j}(o,r)$, $r \geq1$.

This implies $V^{\widetilde M_j}(o,r) \leq V_j(r)$ and $S^{\widetilde M_j}(o,r) \geq \frac12 S_j(r)$ for all $r\geq 1$.
Gluing $\widetilde M_1$ and $\widetilde M_2$ along their boundary we obtain a metric on the connected sum $M_1 \sharp M_2$ such that
\begin{eqnarray*}
\int_1^\infty \frac{V^{M_1 \sharp M_2}(o,r)}{S^{M_1 \sharp M_2}(o,r)}dr
&=&
\int_1^\infty \frac{V^{\widetilde M_1}(o,r) + V^{\widetilde M_2}(o,r)}{S^{\widetilde M_1}(o,r) + S^{\widetilde M_2}(o,r)} dr \\
&\leq&
2\int_1^\infty \frac{V_1(r) + V_2(r)}{S_1(r) + S_2(r)} dr 
\ep<\ep
\infty .
\end{eqnarray*}
Thus we have constructed a 2-dimensional stochastically complete connected manifold such that \eqref{eq:VSendl} holds.
In fact, the manifold has recurrent Brownian motion even.
An easy modification of this construction yields such examples also in dimensions $n \geq3$.
\end{rem}

\end{document}